\documentclass[11pt]{amsart}
\usepackage{amsmath,amssymb,amsthm,amsfonts,enumerate,hyperref, caption,array}
\usepackage{cite}
\theoremstyle{plain}
\newtheorem{propn}{Proposition}[section]
\newtheorem{thm}[propn]{Theorem}
\newtheorem{lemma}[propn]{Lemma}

\theoremstyle{definition}
\newtheorem{defn}[propn]{Definition}

\theoremstyle{remark}

\numberwithin{equation}{section}

\begin{document}

\title{Approximating Ricci solitons and Quasi-Einstein metrics on toric surfaces}
\author{Stuart James Hall}
\address{Department of Applied Computing, University of Buckingham, Hunter St., Buckingham, MK18 1G, U.K.} 
\email{stuart.hall@buckingham.ac.uk}

\author{Thomas Murphy }
\address{Department of Mathematics, California State University Fullerton, 800 N. State College Bld., Fullerton, CA 92831, USA.}
\email{tmurphy@fullerton.edu}

\maketitle

\begin{abstract}
We present a general numerical method for investigating prescribed Ricci curvature problems on toric K\"ahler manifolds. This method is applied to two generalisations of Einstein metrics, namely Ricci solitons and quasi-Einstein metrics. We begin by recovering the Koiso--Cao soliton and the L\"u--Page--Pope quasi-Einstein metrics on $\mathbb{CP}^{2}\sharp\overline{\mathbb{CP}}^{2}$ (in both cases the metrics are known explicitly). We also find numerical approximations to the Wang--Zhu soliton on $\mathbb{CP}^{2}\sharp 2\overline{\mathbb{CP}}^{2}$  (here the metric is not known explicitly). Finally, a substantial numerical investigation of the quasi-Einstein equation  on $\mathbb{CP}^{2}\sharp 2\overline{\mathbb{CP}}^{2}$ is conducted. In this case it is an open problem as to whether such metrics exist on this manifold. We find metrics that solve the quasi-Einstein equation to the same degree of accuracy as the approximations to the Wang--Zhu soliton solve the Ricci soliton equation. 
\end{abstract}
\section{Introduction}
\subsection{Main ideas and results}
 In this paper, the existence of  two related generalisations of Einstein metrics are investigated numerically. Let $(M,g)$ be a complete Riemannian manifold. The metric is called a \emph{gradient Ricci soliton} if  
 \begin{equation}\label{e1}
Ric(g)+\nabla^{2}\phi = \lambda g,
\end{equation}
for a smooth function $\phi \in C^{\infty}(M)$ and constant $\lambda \in \mathbb{R}.$
It is called a \textit{quasi-Einstein metric} if 
\begin{equation}\label{e2}
Ric (g)+\nabla^{2}\phi-\frac{1}{m}d\phi\otimes d\phi = \lambda g,
\end{equation}
for a smooth function $\phi \in C^{\infty}(M)$ and constants $\lambda,m \in \mathbb{R}$ with $m>0$.
If the function $\phi$ is constant then Equations (\ref{e1}) and (\ref{e2}) reduce to that of an Einstein metric. Ricci solitons and quasi-Einstein metrics with nonconstant $\phi$ are referred to as \emph{non-trivial}. Ricci solitons arise as fixed points (up to scaling and diffeomorphism) and as singularity models for the Ricci flow. When $m\in \mathbb{N}$, quasi-Einstein metrics  occur as the base in the warped-product construction of Einstein metrics where the fibre is an Einstein manifold of dimension $m$.\\ 
\\
In dimension 4 the only known compact, non-trivial Ricci solitons and quasi-Einstein metrics occur on the complex surfaces  $\mathbb{CP}^{2}\sharp  \overline{\mathbb{CP}}^{2}$ and $\mathbb{CP}^{2}\sharp 2 \overline{\mathbb{CP}}^{2}$. On  $\mathbb{CP}^{2}\sharp  \overline{\mathbb{CP}}^{2}$ there is a Ricci soliton discovered independently by Koiso \cite{Koi} and Cao \cite{Cao} and, for any $m>1$ a quasi-Einstein metric found by L\"u, Page and Pope \cite{LPP}.  
On $\mathbb{CP}^{2}\sharp 2 \overline{\mathbb{CP}}^{2}$ there is a Ricci soliton found by Wang and Zhu \cite{WZ}. The existence of a quasi-Einstein metric  on $\mathbb{CP}^{2}\sharp 2 \overline{\mathbb{CP}}^{2}$ is an open problem.\\
\\
In Sections \ref{resultsRS} and \ref{resultsQEM}, the results of various algorithms, developed to produce numerical approximations to the  Koiso--Cao, Wang--Zhu and L\"u--Page--Pope metrics, are presented. The more substantial contribution of the paper is to then apply these algorithms to search for quasi-Einstein metrics on $\mathbb{CP}^{2}\sharp 2 \overline{\mathbb{CP}}^{2}$. This is also carried out in Section \ref{resultsQEM}. The extent to which the results suggest the existence of quasi-Einstein metrics on $\mathbb{CP}^{2}\sharp 2 \overline{\mathbb{CP}}^{2}$ is as follows:\vspace{.2in}

\emph{For} $m>1$, \emph{there are conformally K\"ahler metrics on} $\mathbb{CP}^{2}\sharp 2\overline{\mathbb{CP}}^{2}$ \emph{that solve the quasi-Einstein equation (\ref{e2}) with an error that is comparable to the error in solving the Ricci soliton equation (\ref{e1}) by equivalent approximations to the Wang--Zhu Ricci soliton. In particular, there exists an approximately warped-product Einstein metric on }  ${(\mathbb{CP}^{2}\sharp 2\overline{\mathbb{CP}}^{2})\times\mathbb{CP}^{1}}$. \emph{There is numerical evidence that if conformally K\"ahler quasi-Einstein metrics with $J$-invariant Ricci tensor exist on} $\mathbb{CP}^{2}\sharp 2\overline{\mathbb{CP}}^{2}$, \emph{the K\"ahler class of the metric is not the first Chern class.} 

\subsection{Exisiting numerical work}

The Koiso--Cao and L\"u--Page--Pope metrics admit a cohomogeneity one action by $U(2)$ which reduces Equations (\ref{e1}) and (\ref{e2}) to a system of ordinary differential equations. These equations can be solved explicitly  (see \cite{BHJM} for a unified description of these metrics).  Wang and Zhu actually found Ricci solitons on a general class of complex manifolds, namely toric K\"ahler manifolds. Their existence proof is not constructive and so geometers and physicists have been interested in finding numerical approximations of these Ricci solitons. Numerical approximations to the Wang--Zhu soliton were first found by Headrick and Wiseman \cite{HW} by simulating the K\"ahler--Ricci flow which has the Wang--Zhu soliton as a unique fixed point up to automorphisms. Another approximation was found by the first author in his doctoral thesis \cite{HThes} by using the theory of canonically balanced metrics initiated by Donaldson \cite{Donnum}. This fits into a wider program of numerically approximating distinguished Riemannian metrics on K\"ahler manifolds \cite{DonBun}, \cite{DHHJW}, \cite{HWK3} and \cite{Keller}. Non-trivial quasi-Einstein metrics are never K\"ahler \cite{CSW}. However the L\"u--Page--Pope metrics are \emph{conformal} to a toric K\"ahler metric \cite{BHJM} (there is also related work on this topic by Maschler \cite{Mas}).  The methods introduced in this paper are rooted in the theory of toric K\"ahler manifolds, and are similar to those used by the authors in \cite{HM14} to approximate the Chen--LeBrun--Weber Einstein metric on $\mathbb{CP}^{2}\sharp 2 \overline{\mathbb{CP}}^{2}$ \cite{CLW}.  These new algorithms are straightforward to implement in Matlab.   We are confident that they  will generalize to yield numerical approximations to the Wang--Zhu solitons in higher dimensions. 

\subsection{Computer Code}

The computer calculations were carried out using the Matlab optimisation toolbox, in particular, the function `lsqnonlin'.  The various functions used to perform the algorithms are available on the authors' webpages \footnote{\text{http://www.buckingham.ac.uk/directory/dr-stuart-hall/}} \footnote{\text{http://mathfaculty.fullerton.edu/tmurphy/research.html}}. Readers wishing to implement the procedures described should save the functions to a directory where Matlab can access them.  Then, given an initial vector of inputs `x0', call the  `lsqnonlin' function by typing \verb|x = lsqnonlin(@WZT1,x0)| in to the command line.  This will perform the $T_{1}$ minimisation algorithm searching for the Wang--Zhu soliton (see Section \ref{S3} for more details).  Other algorithms can be performed by changing the `WZT1' input to the appropriate function.\\
\\
\textit{Acknowledgements:}\\
The authors are very grateful to the anonymous referee of our related paper \cite{HM14} for suggesting we use the Levenberg-Marquardt algorithm. A good deal of the work was conducted whilst SH was visiting TM at California State University Fullerton. SH would like to thank the department for their hospitality. This trip was funded by a CSUF Startup grant.
\section{Background}

\subsection{Toric K\"ahler metrics}

Key to the numerical algorithms is the theory of toric K\"ahler manifolds. We refer the reader to \cite{A2} and \cite{DonTor} for background. For our purposes, a toric K\"ahler manifold will be a K\"ahler manifold $(M^{2n},\omega, J)$ that admits an effective action of the torus $\mathbb{T}^{n}$ that is simultaneously holomorphic and Hamiltonian. Crucially, there is a dense open subset $M^{o}\subset M$ on which this action is free. From the machinery developed in \cite{A2} and \cite{G} we obtain: 
\begin{itemize}
\item a compact convex polytope $P\subset \mathbb{R}^{n}$(called the moment polytope) and an identification
$$
{M^{o} \cong P^{o}\times \mathbb{T}^{n}},
$$
\item a finite set of  affine linear functions $l_{i}:\mathbb{R}^{n}\rightarrow \mathbb{R}$  such that the polytope $P$ is obtained as the intersection of the regions defined by $l_{i}(x)\geq 0$,
\item a smooth convex function $u:P^{o}\rightarrow \mathbb{R}$ such that in the coordinates $(x_{i},\theta_{j})$ on  $P^{o}\times \mathbb{T}^{n}$ the metric $g(\cdot,\cdot)=\omega(J\cdot,\cdot)$ has the form
$$g = u_{ij}dx_{i}dx_{j}+u^{ij}d\theta_{i} d\theta_{j}.$$
\end{itemize}
Here $u_{ij} = \frac{\partial^{2} u}{\partial x_{i} \partial x_{j}}$ and $u^{ij}$ is the ordinary matrix inverse. In the polytope coordinates the symplectic form $\omega$ is standard and is given by
$$
\omega = \sum_{i}dx_{i} \wedge d\theta_{i}.
$$ 
The metric has coordinate singularities on the boundary $\partial P$ of $P$. It is known exactly how this has to occur. A result of Guillemin \cite{G}  and Abreu \cite{A2} shows that symplectic potential can be written as
\begin{equation}\label{e4}
u = \frac{1}{2}\sum_{i}l_{i}\log(l_{i})+F,
\end{equation}
where $F$ is a smooth function on $P$. The term $\frac{1}{2}\sum_{i}l_{i}\log(l_{i})$ is known as the canonical symplectic potential associated to $P$ and we will denote this by $u_{can}.$ Functions of the form (\ref{e4}) are said to satisfy the \emph{Guillemin boundary conditions}.\\
\\
One nice aspect of toric K\"ahler metrics (especially expounded by Abreu) is that various curvatures have particularly compact expressions in the polytope coordinates.  The Ricci curvature in the polytope coordinates can be computed as
\begin{equation}\label{Ricci}
Ric(\partial_{x_{i}},\partial_{x_{j}}) = \frac{1}{2}\left( \frac{\partial^{2}}{\partial x_{i} \partial x_{j}} -u^{kl}\frac{\partial u_{ij}}{\partial x_{k}}\frac{\partial}{\partial x_{k}}\right)\log( \det (D^{2} u))
\end{equation}
where $D^{2}u$ is the Euclidean Hessian of $u$. This completely determines the Ricci curvature as the `mixed' terms ${Ric(\partial_{x_{i}},\partial_{\theta_{j}})=0}$ and other terms can be computed using $J$-invariance.  The scalar curvature is given by
\begin{equation}
S = -\sum_{i,j}u^{ij}_{ij}.
\end{equation}

The polytope associated to toric K\"ahler metrics on $\mathbb{CP}^2\sharp\overline{\mathbb{CP}}^{2}$ is the trapezium determined by the affine linear functions 
$$l_{1}(x) = a+x_{1}+x_{2}, \ l_{2}(x) = 1+x_{1}, \ l_{3}(x) =1+x_{2} \textrm{ and } l_{4}(x)=1-x_1-x_2.$$ 
The cohomology class of the metric is determined by the parameter \newline $a \in (-1,2)$.\\
\\
The toric K\"ahler metrics we will be interested in on $\mathbb{CP}^{2} \sharp2\overline{\mathbb{CP}}^{2}$ have associated polytope given by
\begin{align*}
l_{1}(x) =& 1+x_{1}, \ l_{2}(x) = 1+x_{2}, \ l_{3} = a-1 -x_{1}, \ l_{4} = a-1-x_{2}, \text{ and}\\
l_{5}(x) =& a-1-x_{1}-x_{2}.
\end{align*}
Again, the cohomology class of the metric is determined by the parameter $a$.
\subsection{K\"ahler--Ricci Solitons}
If a gradient Ricci soliton is  also a K\"ahler metric then it is referred to as a \textit{K\"ahler--Ricci soliton} and there are a number of properties that can be deduced from Equation (\ref{e1}).  Since the Ricci tensor and the metric are  both $J$-invariant, so is the Hessian $\nabla^{2}\phi$. This implies that $\nabla \phi$ is a real holomorphic vector field. Scaling so that $\lambda=1$ in Equation (\ref{e1}) and using the complex structure,  the equation for the associated  $(1,1)$-forms  is
$$\rho+i\partial\bar{\partial} \phi = \omega,$$
where $\rho$ is the Ricci form. Hence $[\omega]=2\pi c_{1}$ and so the manifold must be a smooth Fano variety. The Koiso--Cao and Wang--Zhu solitons are toric K\"ahler metrics. Their existence is a special case of; 
\begin{thm}[Wang--Zhu \cite{WZ}]
Let $(M,J)$ be a Fano toric K\"ahler manifold. Then there exists a K\"ahler-Ricci soliton on $M$ unique up to automorphism.
\end{thm}
In the equations defining the polytopes,  the Koiso--Cao metric has $a=1$ and the Wang--Zhu metric has $a=2$. \\
\\
The potential function $\phi$ of a toric K\"ahler Ricci soliton can be taken to be $\mathbb{T}^{n}$-invariant and also a function of the polytope coordinates $x_{1},..., x_{n}$. Calculation shows that if the Hessian $\nabla^{2}\phi$ is $J$-invariant then $\phi$ must be an affine linear function in the polytope coordinates;
$$\phi(x_{1},...,x_{n}) = \sum_{i=1}^{i=n}a_{i}x_{i}.$$
Using methods contained in  \cite{DonTor}, \cite{WZ}, it is possible to work out the coefficients $a_i$ without explicit knowledge of the metric. This was done for the Koiso--Cao and Wang--Zhu metrics in \cite{Halldga}.  Specifically, 
$$
\phi_{KC} = 0.527620(x_1+x_2) \text{ and } \phi_{WZ} = -0.434748(x_1+x_2)
$$
where $\phi_{KC}$ and $\phi_{WZ}$ denotes the potential for the Koiso--Cao metric and Wang--Zhu metrics respectively and the coefficents are given to 6 significant figures. \\ 
\\
As mentioned in the introduction, the Koiso--Cao metric admits an isometric $U(2)$ action. This means the function $F$ can be taken as a function of $t = x_1+x_2$. The second derivative of the function $F$ (as this is all that is needed to determine the metric) was calculated in \cite{BHJM} to be
\begin{equation}\label{e5}
F''(t) = \left(\frac{\frac{1}{2}c^{3}(2+t)}{c^{3}de^{c(2+t)}+c^{2}t(2+t)+2c(1+t)+2} +\frac{\frac{1}{2}(t^{2}-2t-5)}{(1-t^{2})(t+2)}\right)
\end{equation}
where $c\approx 0.527620$ is the coefficient in the potential function and ${d\approx -6.91561}$.

The algorithms approximating Ricci solitons exploit the following fact:
\begin{lemma}\label{trrseqn}
Let $(M, \omega, J)$ be a K\"ahler metric of real dimension $n$, such that $\omega \in 2\pi c_1(M)$. If $\phi\in C^{\infty}(M)$ has holomorphic gradient  and solves
\begin{equation}\label{scalrsol}
S + \Delta \phi = n
\end{equation}
then $(M,\omega, \phi)$ is a K\"ahler-Ricci soliton. 
\end{lemma}
\proof By assumption, there exists a function $\eta$ such that $$\omega = \rho + i\partial\overline{\partial}\eta.$$
Taking the trace, we obtain $$ n = S + \Delta\eta.$$ The result follows from the uniqueness (up to a constant) of solutions to Poisson's equation. 
\endproof 

\subsection{Quasi-Einstein metrics}   

The quasi-Einstein equation has a first integral due to Kim--Kim \cite{KK} coming from the contracted second Bianchi identity. For any solution to Equation (\ref{e2}), there is a constant $\mu$ for which the quasi-Einstein potential $\phi$ satisfies
\begin{equation}\label{kkeqn}
1-\frac{1}{m}\left( \Delta\phi - |\phi|^2\right)= \mu e^{\frac{2\phi}{m}}.
\end{equation}
As mentioned in the introduction, this (together with Equation (\ref{e2}))  implies that when $m\in \mathbb{N}$ and $(F^m,h)$ is an Einstein metric with Einstein constant $\mu$, ${(M\times F, g \oplus e^{\frac{-2f}{m}}h)}$ is an Einstein metric with Einstein constant $\lambda$.\\ 
\\
At the time of writing, the only general construction of compact quasi-Einstein metrics where ${m >1}$ is an arbitrary real number are the L\"u--Page--Pope examples (and generalizations due to the first author \cite{Halljgp}).   A foundational result due to Case, Shu and Wei \cite{CSW} states that non-trivial quasi-Einstein metrics are never K\"ahler. However, the L\"u--Page--Pope metrics $g_{LPP}$ are \emph{conformal} to toric K\"ahler metrics, and the K\"ahler metrics lie in the first Chern class \cite{BHJM} (so $a=1$). We will write $g_{LPP} = e^{2\sigma}g_K$, where $g_K$ is the toric K\"ahler metric.   In \cite{BHJM}, an explicit description of  $g_{LPP}$ is give. This is possible due to the fact the metrics are invariant under a cohomogeneity one action of  $U(2)$. In particular, the $U(2)$ invariance forces $g_{LPP}$ to have $J$-invariant Ricci tensor. All this implies that the conformal function $\sigma$ and the quasi-Einstein potential $\phi$, viewed as functions on the trapezium, are given as
\begin{equation}\label{e3}
\sigma = -\log (bt + c) \text{ and } \phi = -m \log (\frac{dbt +dc +1}{bt+c}),
\end{equation}
where $b,c$ and $d$ are constants and $t=x_1+x_2$. Given $m$, the constants $b,c,d$ can be determined by
$$
d = (2(2b - c))^{-1} \text{ and } c^2 = b^2+1.
$$
There is an additional constraint;
$$\int_{-1}^{1} \frac{(2+s) - 2(bs+c)^2(dbs + dc+ 1)^{m-2}(2+s)}{(bs+c)^{m+4}}ds = 0 
$$
For example, when $m=2$, we find $b \approx 0.076527$, $c\approx 1.002924$ and $d\approx 0.588325$.\\
\\
It is natural to begin searching for quasi-Einstein metrics on  $\mathbb{CP}^2\sharp2\overline{\mathbb{CP}}^2$ by looking for metrics with similar properties to the L\"u--Page--Pope examples. For example, looking for quasi-Einstein metrics that are conformal to toric K\"ahler metrics and that have $J$-invariant Ricci tensor. If such metrics were to exist,  the conformal function $\sigma$ and quasi-Einstein potential $\phi$ would also have the form given by Equation (\ref{e3}). In \cite{BHJM}, a set of constraints for the parameters $a,b,c,d$ and $\mu$ (the constant appearing in Equation (\ref{kkeqn})) was also derived: namely
$$ \frac{4b}{(c-2b)(dc+1-2db)} =\frac{1}{(c-2b)^2}-\frac{\mu}{(dc+1-2db)^2}, $$
$$0 = \frac{1}{(c+(a-2)b)^2}-\frac{\mu}{(dc+1+(a-2)db)^2}, $$
and
$$\frac{-2b}{(c+(a-1)b)(dc+1+(a-1)db)}=\frac{1}{(c+(a-1)b)^{2}}-\frac{\mu}{(dc+1+(a-1)db)^{2}} .$$
Moreover we have the following;
$$\int_{P} (e^{-\phi}-\mu e^{\left(\frac{2}{m}-1\right)\phi})e^{4\sigma}dx = 0. $$
If the K\"ahler metric were in the first Chern class  then $a=2$. The constraint equations can then be numerically solved for each fixed $m$. For example, when $m=2$, $b\approx -0.0744357$, $c\approx 1.00482$, $d\approx -0.463585$, and $\mu \approx 0.282687$.\\ 
\\
A useful quantity that will be utilised in the algorithm is the conformal quasi-Einstein equation that is satisfied by the K\"ahler metric $g_{K}$, the conformal function $\sigma$ and the potential function $\phi$:
\begin{equation}\label{confqem}
 Ric(g_K) = \mathcal{A}(g_K, \phi, \sigma, m)
\end{equation}
where the right hand side is the tensor
\begin{align*}
 \mathcal{A}(g_K, \phi, \sigma, m) =  &-2\nabla^2\sigma + \nabla^2\phi + 4d\sigma\otimes d\sigma - d\sigma\otimes d\phi - d\phi\otimes d\sigma \\ &-\frac{1}{m} d\phi\otimes d\phi + (g_K(\nabla\sigma, \nabla\phi) - 2|\nabla\sigma|^2 - \Delta\sigma - e^{2\sigma})g_K.
\end{align*}
where all the geometric quantities are calculated with respect to the metric $g_{K}$.

\section{The key ideas of the approximation algorithm}\label{S3}

The numerical approximations of the metrics are presented as symplectic potentials. The following space was essentially introduced by Doran et. al. \cite{DHHJW} to give approximations to Siu's K\"ahler-Einstein metric on $\mathbb{CP}^{2} \sharp 3\overline{\mathbb{CP}}^{2}$. 
  
\begin{defn} The space of \emph{restricted symplectic metrics} of degree $d$, $\mathcal{S}_d$ is  defined as the  metrics given by  a symplectic potential where the function $F$ in Equation (\ref{e4}) is a polynomial of degree $d$ in the polytope coordinates. 
\end{defn}

It is clear that $\mathcal{S}_{d}$ is an open subset of $\mathbb{R}^{(d+1)^{2}}$. In fact, the metrics we are interested in are also invariant under a $\mathbb{Z}_{2}$ action switching $x_{1}$ and $x_{2}$ and so we will work with the set of $\mathbb{Z}_{2}$-invariant restricted symplectic metrics of degree $d$, $\mathcal{S}_{d}^{\mathbb{Z}_{2}}$. Practically, this means that the function $F$ is determined by the ${N_{d} = \lfloor \frac{d^{2}+6d+1}{4}\rfloor}$ coefficients $c_{1},...,c_{N_{d}}$ by taking
\begin{equation}\label{Ftaylor}
F(x_{1},x_{2}) = c_{1}x_{1}x_{2}+c_{2}(x_{1}^{2}+x_{2}^{2})+...+c_{N_{d}}(x_{1}^{d}+x_{2}^{d}).
\end{equation}
All of the algorithms involve minimizing a function 
$$
\mathfrak{F}: \mathcal{S}_{d}^{\mathbb{Z}_{2}}\rightarrow \mathbb{R},
$$ 
given by
\begin{equation}\label{FrakFeq}
\mathfrak{F}(g) = \int_P (T(g))^{2} dx_1 \ dx_2
\end{equation}
where $T(g)$ is a geometric quantity that depends upon the restricted symplectic metric. We use algorithms that attempt to minimise a scalar curvature type quantity and a Ricci curvature type quantity.\\ 
\\
The Ricci curvature type integrands $T(g)$ require the evalutation of quantities (such as the coeffcients of the metric $g$) that become singular at the boundary $\partial P$ of the polytope $P$. In this case the integrals are computed over the polytope with parallel boundary $P_{\delta}$ defined by taking ${l_{r}(x_{1},x_{2})>\delta>0}$ where $l_{r}$ are the affine linear functions defining $P$. We have implicitly assumed that the symplectic potentials of K\"ahler Ricci solitons and conformally K\"ahler quasi-Einstein metrics are analytic in the polytope coordinates (a justification in the case of Ricci solitons is given in \cite{DHHJW}).  Hence if one can solve Equations (\ref{e1}) and (\ref{e2}) on an open set of the polytope this should completely characterise the metric.  For this reason the choice of $\delta$ is not too important as it only controls numerical error. We found taking $\delta=0.005$ yielded good results.  
\subsection{Approximating Ricci solitons}

The algorithms used to approximate Ricci solitons employ two choices of $T$. Firstly, one can choose to minimize the functional given by taking 
$$
T_1(g)= S+\Delta \phi-4,
$$
in Equation (\ref{FrakFeq}).  As explained in Lemma \ref{trrseqn}, if a K\"ahler metric satisfies $T_1(g) \equiv 0$ then it also solves the soliton equation (\ref{e1}). Secondly one could choose
$$
T_2(g) := \sqrt{\sum_{1\leq i,j \leq 2}(Ric_{ij}+\nabla^{2}\phi_{ij}-u_{ij})^{2}},
$$
in Equation (\ref{FrakFeq}). The terms of the sum in $T_2(g)$ become singular at the boundary $\partial P$.  In this case the integral (and hence the least squares function) uses the polytope $P_{0.005}$.

\subsection{Searching for Quasi-Einstein metrics}
To search for quasi-Einstein metrics, the algorithms minimised the functional given by taking the Ricci curvature-type quantity (c.f. Equation (\ref{confqem}))
$$T_3(g) := \sqrt{\sum_{1\leq i,j \leq 2}(Ric_{ij}(g_K)  - \mathcal{A}(g_K, \phi, \sigma, m)_{ij})^2},$$
in Equation (\ref{FrakFeq}). As in the case with the Ricci soliton algorithm, the polytope used in the integral is $P_{0.005}$. If $T_3(g)\equiv 0$ then $e^{2\sigma}g_{K}$ must be a quasi-Einstein metric.\\
\\
There is a related scalar curvature quantity given by taking the trace of Equation (\ref{e2}) to obtain
$$T_4(g): = S + \Delta\phi -\frac{1}{m}|\nabla\phi|^2 - 4.$$
However it is not known (and possibly not true) that a metric solving \newline $T_4(g)\equiv 0$ must be a solution of Equation (\ref{e2}). Hence in the case of searching for quasi-Einstein metrics, the integral of this quantity is used simply as an indication of the accuracy of the approximation.

\subsection{Approximating the integrals}

The integrals were approximated using Gaussian quadrature. For a one-dimensional integral (normalised so that the range is $[-1,1]$) the idea of Gaussian quadrature is to approximate the integral by taking a weighted sum of values
$$\int_{-1}^{1} f(t) dt \approx \sum_{i=1}^{i=k}w_{i}f(t_{i}).$$
The points $t_{i}$ at which the function are sampled are known as the abscissa and the $w_{i}$ are refered to the weights. The points $t_{i}$ and weights $w_{i}$ are chosen so that if $f$ is a polynomial of degree $2k-1$ or less, then the sum will compute the integral exactly.\\ 
\\
To compute integrals over the trapezium $P$ corresponding to the toric metrics on $\mathbb{CP}^{2}\sharp  \overline{\mathbb{CP}}^{2}$ the splitting
$$
\int_P f dx_1 \ dx_2 = \int_{-1}^{1+a}\int_{a-x_1}^{1-x_1} f dx_2 \ dx_1 + \int_{1+a}^{2}\int_{-1}^{1-x_1} f dx_2 \ dx_1,
$$
was taken. Similarly for the pentagon $P$ corresponding to the toric metric on $\mathbb{CP}^{2}\sharp  2\overline{\mathbb{CP}}^{2}$ the following splitting:
$$\int_{P}fdx_{1} \ dx_{2} = \int_{-1}^{a-1}\int_{-1}^{1}fdx_{1} \ dx_{2}+\int_{1}^{a-1}\int_{-1}^{a-1-x_{1}}fdx_{2} \ dx_{1},$$
was taken. These splittings were taken to ensure that all the functions in the one-dimensional interated integrals were smooth. Similar splittings were used when integrating over the polytopes $P_{0.005}$ The iterated one-dimensional integrals were approximated using the Gaussian quadrature method. We took 20 points in the one-dimensional integrals.\\

\subsection{Nonlinear least squares problems}

What the Gaussian quadrature method amounts to is the approximation of the functions $\mathfrak{F}$ by a sum of squares, 
$$\mathfrak{F}(g) \approx \sum_{i=1}^{i=800}(\tilde{w}_{i}T(g))^{2}(p_{i},q_{i}),$$
where $(p_{i},q_{i}) \in P$ are the points of the polytope used in the quadrature procedure and $\tilde{w}_{i}$ is a weight coming from the abscissa weights and the transformation of the one-dimensional integrals. Hence a good approximation to the minimum of $\mathfrak{F}$ over the set $\mathcal{S}^{\mathbb{Z}_{2}}_{d}$ can be found by minimising the function
$$\mathcal{I}(c_{1},...,c_{N_{d}}) = \sum_{i=1}^{i=800}(\tilde{w}_{i}T(c_{1},...,c_{N_{d}}))^{2}(p_{i},q_{i}).$$
This fits in to the framework of \emph{nonlinear least squares problems}. The optimisation toolbox in Matlab has a variety of inbuilt methods for finding approximate solutions of such problems. We used the Levenberg-Marquardt algorithm which we give an overview of here. The problem is to minimise a function
\begin{equation}\label{residual}
\chi(c) = \sum_{i=1}^{i=N}(y_{i}-y_{i}(c))^{2}
\end{equation}
where $c\in \mathbb{R}^{m}$ and $y_{i}:\mathbb{R}^{m}\rightarrow \mathbb{R}$ are functions of $c$. If, near a minimum point $c^{\ast}$, the function $\chi$ is well approximated by its quadratic Taylor expansion, an initial guess of the minima $c_{in}$ can be updated via
\begin{equation}\label{U1}
c_{new} = c_{in}+\left(\nabla^{2}\chi (c_{in})\right)^{-1}(\nabla \chi(c_{in})).
\end{equation}
If this is not the case then the initial guess can be updated via the gradient descent method
\begin{equation}\label{U2}
c_{new} = c_{in} -\gamma (\nabla \chi(c_{in})),
\end{equation}
for an appropriate constant $\gamma >0$ which is determined by the Hessian $\nabla \chi(c_{in})$ . The fact that $\chi$ is a sum of squares allows an approximation of the Hessian close to $c^{\ast}$
$$(\nabla^{2}\chi)_{kl} \approx 2\sum_{i=1}^{i=N}\frac{\partial y_{i}}{\partial c_{k}}\frac{\partial y_{i}}{\partial c_{l}}.
$$
The Levenberg-Marquardt algorithm compares the residual (\ref{residual}) at the initial point $c_{in}$ and at the update $c_{new}$ given by (\ref{U1}). If there is no improvement, then one step in the gradient descent is performed with (\ref{U2}).  This process is explained in more detail in Section 15.5.2 of \cite{NM}. It is worth noting that this algorithm does have the problem that it can become trapped in `flat' regions near to to the minimum or even in local minima. Other methods may need to be used in conjunction in order to find global minima.\\
\\
The Levenberg-Marquardt method for minimising $\mathcal{I}$ was implemented in Matlab.  To do this an initial guess for the coefficients $(c_{1},...,c_{N_{d}})$ and conditions for the algorithm to terminate the search needed to be specified. Usually a search was started in the space of quadratic approximations $\mathcal{S}^{\mathbb{Z}_{2}}_{2}$ with $F=0$. Subsequent searches in $\mathcal{S}^{\mathbb{Z}_{2}}_{d}$ took the approximation found in $\mathcal{S}^{\mathbb{Z}_{2}}_{d-1}$ as the inital condition. In general, the following stopping conditions were used for the algorithm:
\begin{enumerate}
\item The algorithm stops when there have been 4000 function evaluations.
\item The algorithm stops if the absolute value of the change in the residual is less than ${5\times 10^{-12}}$.
\item The algorithm stops if the Euclidean norm of the change in the vector $c$ is less than ${5\times 10^{-12}}$.
\end{enumerate}
Where the search procedure differed from this, the specific method is explained with the results.

\subsection{Error estimates}
To determine the accuracy of the approximation given by the vector ${(c_{1}^{\ast},...,c_{N_{d}}^{\ast})}$ returned by the minimisation procedure, various measures of how far the numerical approximation fails to solve Equations (\ref{e1}) or (\ref{e2}) were computed. 

\begin{defn}\label{normalizederror}
The normalized error associated to $T$ at $(c_{1}^{\ast},...,c_{N_{d}}^{\ast})$  is defined as
$$
\mathcal{E}\left(T(c_1^{\ast},...., c_{N_d}^{\ast})\right) := Vol(P)^{-1}\sqrt{\mathcal{I}(c_1^{\ast},...., c_{N_d}^{\ast})}.
$$
The absolute error associated to $T$ at $(c_{1}^{\ast},...,c_{N_{d}}^{\ast})$ is
$$
Max\left(T(c_{1}^{\ast},...,c_{N_{d}}^{\ast})\right) := \max_{P_{\delta}}|T(c_{1}^{\ast},...,c_{N_{d}}^{\ast})|
$$

\end{defn}

The quantities $\mathcal{E}$  are essentially the residual errors for each of the least square problems.

\subsection{Further refinements to the algorithm}

There is no \emph{a priori} reason for a toric K\"ahler metric  on  $\mathbb{CP}^2\sharp 2\overline{\mathbb{CP}}^2$, conformal to a quasi-Einstein metric, to  be in the first Chern class.  Hence when searching for quasi-Einstein metrics using the preceding algorithms, it is necessary to view the parameter $a$ (and hence $b$, $c$ and $d$) as unknown variables to be determined. For the Koiso--Cao and Lu--Page--Pope metrics, where these variables are \emph{a priori} known, one can still run an unconstrained search and check that the algorithm finds the correct values. \\
\\
Secondly, there is also no reason to suppose that a conformally K\"ahler metric on $\mathbb{CP}^2\sharp 2\overline{\mathbb{CP}}^2$ must have $J$-invariant Ricci tensor. However if there is a continuous family parameterised by $m$ and the limit of such a family as $m$ tends to infinity is the Wang-Zhu soliton, then, for large values of $m$ at least, the $J$-anti-invariant part of the Ricci tensor cannot be too large. For this reason we also consider searching for quasi-Einstein metrics $g=e^{2\sigma}g_{K}$ where
$$
\sigma = -\log(bt+c) +\epsilon_{1}(x_{1},x_{2}) \textrm{ and } \phi = -m\log\left(\frac{dbt +dc+1}{bt+c}\right)+\epsilon_{2}(x_{1},x_{2}).
$$
The functions $\epsilon_{i}(x_{1},x_{2})$ can both be expanded as a $\mathbb{Z}_{2}$-invariant polynomial in the $x_{1}$ and $x_{2}$ coordinates.
\section{results: Ricci solitons} \label{resultsRS}

\subsection{The Koiso--Cao soliton}
As the symplectic potential for the Koiso--Cao soliton is known explicitly, it gives an important first check as to the accuracy and reliability of the proposed algorithms. Recall that Equation (\ref{e5}) gives $F$ explicitly as a function of the quantity $t=x_{1}+x_{2}$.  It is possible to expand $F$ as a Taylor series in $t$,
$$F(t) = k_{1}t^{2}+k_{2}t^{3}+...$$
Using the Equation (\ref{e5}), Table \ref{T1} gives the coefficients $k_{1},...,k_{4}$ to 6 significant figures.
\begin{table}[!ht] 
\centering
\caption{The first four coefficients of the Taylor expansion of $F$}
\begin{tabular}{|c|c|c|c|c|}
\hline
$k_{1}$ & $k_{2}$ & $k_{3}$ & $k_{4}$ \\ 
\hline
-0.0900384 & 0.0159081 & $-4.25806 \times 10^{-3}$ & $1.34121\times 10^{-3}$\\
\hline
\end{tabular}
\label{T1}
\end{table}

By imposing additional relations on the coefficients $c_{i}$ in Equation (\ref{Ftaylor}) (${c_{1} = 2c_{2}}$, ${c_{3}=3c_{4}}$, etc.), the minimisation procedures can take place in the space of $U(2)$ invariant symplectic potentials $\mathcal{S}^{U(2)}_{d}\subset \mathcal{S}^{\mathbb{Z}_{2}}_{d}$. The $T_{1}$-minimisation algorithm was run on these spaces and the corresponding Taylor coefficients are recorded in Table \ref{T2}. Here the first 4 coefficients of the approximation found in $ \mathcal{S}^{U(2)}_{10}$ are given. There is a good approximation with the coefficients being accurate to five or six decimal places.

\begin{table}[!ht]
\centering
\caption{The first four Taylor coefficients of $F$ for the $S_{10}^{U(2)}$ approximation by the $T_{1}$ method}
\begin{tabular}{|c|c|c|c|}
\hline
$k_{1}$ & $k_{2}$ & $k_{3}$ & $k_{4}$\\ 
\hline
-0.0900413 & 0.0159070 & $-4.25899 \times 10^{-3}$ & $1.34848\times 10^{-3}$\\
\hline
\end{tabular}
\label{T2}
\end{table}

The $T_{1}$-minimisation algorithm was performed on the larger space $\mathcal{S}^{\mathbb{Z}_{2}}_{d}$.  The results of the approximation at each $d$ are given in Table \ref{T3}. Here there is exponential convergence of the error terms toward 0.\\
\\
\begin{table}[!ht]
\centering
\caption{$T_{1}$ minimisation on the spaces $\mathcal{S}^{\mathbb{Z}_{2}}_{d}$ for the Koiso-Cao soliton}
\begin{tabular}{|c|c||c|c||c|c|}
\hline 
Degree $d$ & $N_{d}$ & $\mathcal{E}(T_{1})$  & $Max(T_{1})$ & $\mathcal{E}(T_{2})$ & $Max(T_{2})$  \\
\hline
\hline
2 & 2 & 0.27 &  1.8 &  5.0  & $2.2 \times 10^{2}$ \\
\hline
3 & 4 & 0.12 & 0.82 & 2.3 & 33\\
\hline
4 & 7 & 0.047 & 0.34  & 0.87 & 3.0\\
\hline
5 & 10 & 0.017 & 0.13 & 0.31 & 0.17\\
\hline
6 & 14 & $6.1\times 10^{-3}$ &  0.049 & 0.10 & $8.2\times 10^{-3}$ \\
\hline 
7 & 18 & $2.1\times 10^{-3}$ &  0.017 & 0.034 & $1.7\times 10^{-3}$ \\ 
\hline
8 & 23 & $6.9\times 10^{-4}$ & $6.0\times 10^{-3}$ &  0.010 & $2.2\times 10^{-4}$ \\
\hline
9 & 28 & $2.2 \times 10^{-4}$ &  $2.0\times 10^{-3}$ &  $3.0 \times 10^{-3}$ & $5.1\times 10^{-5}$ \\
\hline
10 & 34 & $7.2 \times 10^{-5}$  &  $6.8 \times 10^{-4}$ & $8.5 \times 10^{-4}$ & $8.7 \times 10^{-6}$\\
\hline
\end{tabular}

\label{T3}
\end{table}

The $T_{2}$-minimisation algorithm was also performed on the spaces $\mathcal{S}^{U(2)}_{d}$.  The first four Taylor coefficients of the approximation found in $\mathcal{S}^{U(2)}_{10}$  are given in Table \ref{T4}. As with the $T_{1}$ algorithm, there is good approximation to the exact coefficients in Table \ref{T1} in this case.\\
\\
\begin{table}[!ht]
\centering
\caption{The first four Taylor coefficients of $F$ for the $\mathcal{S}_{10}^{U(2)}$ approximation by the $T_{2}$ method}
\begin{tabular}{|c|c|c|c|}
\hline
$k_{1}$ & $k_{2}$ & $k_{3}$ & $k_{4}$ \\ 
\hline
 -0.0900268 & 0.0159100 & $-4.29714 \times 10^{-3}$ & $1.35386\times 10^{-3}$ \\
\hline
\end{tabular}
\label{T4}
\end{table}

The $T_{2}$ minimisation method was performed on the larger spaces $\mathcal{S}^{\mathbb{Z}_{2}}_{d}$. The results are contained in Table \ref{T5}.  Here the initial convergence towards 0 is exponential but the algorithm ceases to find a substantial improvement after degree 7.  This is the first indication of the limitation of the method, it seems likely that the algorithm stalls as it is trapped by regions containing local minima.  For comparison, the error profile of the $T_{2}$-minimisation procedure on the spaces $\mathcal{S}^{U(2)}_{d}$ is given in Table \ref{T6}; here the error converges exponentially over the range of degrees considered.

\begin{table}[!ht]
\centering
\caption{$T_{2}$ minimisation on the spaces $\mathcal{S}^{\mathbb{Z}_{2}}_{d}$ for the Koiso--Cao soliton}
\begin{tabular}{|c|c||c|c||c|c|}
\hline 
Degree $d$ & $N_{d}$ & $\mathcal{E}(T_{1})$  & $Max(T_{1})$ & $\mathcal{E}(T_{2})$ & $Max(T_{2})$  \\
\hline
\hline
2 & 2 & 0.33 & 1.4 & 5.2  & $1.2\times 10^{2}$ \\
\hline
3 & 4 & 0.23 & 0.43  & 0.54 & 3.1 \\
\hline
4 & 7 & 0.084  & 0.14  & 0.22 & 0.48 \\
\hline
5 & 10 & 0.044 & 0.082 & 0.072 & 0.038 \\
\hline
6 & 14 & 0.022 & 0.045  & 0.036 & 0.013\\
\hline 
7 & 18 & 0.015 & 0.035 & 0.024 & $6.4 \times 10^{-3}$ \\
\hline
8 & 23 & 0.013 & 0.030 & 0.020  & $3.5 \times 10^{-3}$ \\
\hline
9 & 28 & 0.011  & 0.029 &  0.016 & $2.5 \times 10^{-3}$\\
\hline
10 & 34 & 0.010 &  0.028 & 0.016 & $2.4 \times 10^{-3}$ \\
\hline
\end{tabular}
 
\label{T5}

 \end{table}

\begin{table}[!ht]
\centering
\caption{$T_{2}$ minimisation on the spaces $\mathcal{S}^{U(2)}_{d}$ for the Koiso--Cao soliton}
\begin{tabular}{|c|c||c|c||c|c|}
\hline 
Degree $d$ & $N_{d}$ & $\mathcal{E}(T_{1})$  & $Max(T_{1})$ & $\mathcal{E}(T_{2})$ & $Max(T_{2})$  \\
\hline
\hline
2 & 1 & 0.32 & 1.4  &  5.2 & 120 \\
\hline
3 & 2 & 0.23 & 0.42  & 0.54 & 3.2 \\
\hline
4 & 3 & 0.084 & 0.14 & 0.22 & 0.48 \\
\hline
5 & 4 & 0.042 & 0.076  & 0.072 & 0.043 \\
\hline
6 & 5& 0.015 & 0.028 & 0.025 & $4.5 \times 10^{-3}$ \\
\hline 
7 & 6 & $5.1 \times 10^{-3}$ & $9.4 \times 10^{-3}$ & $8.3\times 10^{-3}$ & $4.3\times 10^{-4}$ \\
\hline
8 & 7 & $1.7 \times 10^{-3}$ & $3.1\times 10^{-3}$ & $2.7 \times 10{-3}$  & $4.0\times 10^{-5}$  \\
\hline
9 & 8 & $6.1\times 10^{-4}$  & $1.2\times 10^{-3}$ &  $9.6\times 10^{-4}$ & $7.6\times 10^{-6}$\\
\hline
10 & 9 & $4.8\times 10^{-4}$ & $9.7\times 10^{-4}$  & $8.8\times 10^{-4}$ & $4.1 \times 10^{-6}$ \\
\hline
\end{tabular}

\label{T6}

\end{table}

\subsection{The Wang--Zhu soliton}

The $T_{1}$-minimisation algorithm was performed on the spaces $\mathcal{S}^{\mathbb{Z}_{2}}_{d}$. The results of this are given in Table \ref{T7}. Here the error term $\mathcal{E}(T_{1})$ seems to be converging exponentially to 0 over the range of degrees considered.  We also give the quartic approximation found in $\mathcal{S}^{\mathbb{Z}_{2}}_{4}$: 
$$
u_{WZ}(x_{1},x_{2}) \approx u_{can}-0.083x_{1}x_{2}-0.121(x_{1}^{2}+x_{2}^{2})-0.038x_{1}x_{2}(x_{1}+x_{2})
$$
$$
-0.029(x_{1}^{3}+x_{2}^{3}) - 0.013x_{1}^{2}x_{2}^{2}-0.010x_{1}x_{2}(x_{1}^{2}+x_{2}^{2})-0.007(x_{1}^{4}+x_{2}^{4}).
$$
We see that the coefficients are broadly similar to the ones found by Headrick and Wiseman in \cite{HW}.  Hence there is good evidence that the succesive approximations are converging to the Wang--Zhu soliton.\\ 
 
\begin{table}[!ht]
\centering
\caption{$T_{1}$  minimisation on the spaces $\mathcal{S}^{\mathbb{Z}_{2}}_{d}$ for the Wang--Zhu soliton}
\begin{tabular}{|c|c||c|c||c|c|}
\hline 
Degree $d$ & $N_{d}$ & $\mathcal{E}(T_{1})$  & $Max(T_{1})$ & $\mathcal{E}(T_{2})$ & $Max(T_{2})$  \\
\hline
\hline
2 & 2 & 0.46 & 3.9 & 6.0  & $0.7\times 10^{3}$ \\
\hline
3 & 4 & 0.24 & 2.5  & 3.0 & $9.6\times 10^{2}$ \\
\hline
4 & 7 & 0.12 &  1.6 & 1.4 & $3.3 \times 10^{2}$  \\
\hline
5 & 10 & 0.064  & 0.90  & 0.70 & 66 \\
\hline
6 & 14 & 0.035 & 0.57 & 0.37 & 21\\
\hline 
7 & 18 & 0.019 & 0.35 & 0.20 & 6.2 \\
\hline
8 & 23 & 0.010 & 0.21  & 0.10  & 1.6 \\
\hline
9 & 28 & $5.5\times 10^{-3}$ & 0.12  &  0.054 & 0.36\\
\hline
10 & 34 & $2.9\times 10^{-3}$  & 0.069 & 0.027 & 0.067\\
\hline
\end{tabular}

\label{T7}
\end{table}
 
Table \ref{T8} contains the results of the $T_{2}$ minimisation procedure on the spaces $\mathcal{S}^{\mathbb{Z}_{2}}_{d}$. As with the $T_{1}$ algorithm, the errors are initially converging exponentially to 0. The quartic approximation obtained here is given by:
$$
u_{WZ} \approx u_{can} - 0.0008x_{1}x_{2}-0.071(x_{1}^{2}+x_{2}^{2})-0.087x_{1}x_{2}(x_{1}+x_{2})-0.048(x_{1}^{3}+x_{2}^{3})
$$
$$
-0.033x_{1}^{2}x_{2}^{2}-0.032x_{1}x_{2}(x_{1}^{2}+x_{2}^{2})-0.022(x_{1}^{4}+x_{2}^{4}).
$$
The values here are somehwat different to those obtained in the $T_{1}$ approximation.  This is probably due to the major contribution to the error being the singular behaviour of the Ricci tensor on the boundary.\\
\\
It will be useful to compare the convergence of this algorithm in higher degrees with that of the algorithms searching for quasi-Einstein metrics on this manifold.  Hence we show the result of continuing the search to the space $\mathcal{S}_{15}^{\ast}$ This represents our best approximation to the Wang--Zhu soliton and the coefficients are on our websites in the file `WangZhu70.txt'.   

 \begin{table}[!ht]
\centering
\caption{$T_{2}$  minimisation on the spaces $\mathcal{S}^{\mathbb{Z}_{2}}_{d}$ for the Wang--Zhu soliton}
\begin{tabular}{|c|c||c|c||c|c|}
\hline 
Degree $d$ & $N_{d}$ & $\mathcal{E}(T_{1})$  & $Max(T_{1})$ & $\mathcal{E}(T_{2})$ & $Max(T_{2})$  \\
\hline
\hline
2 & 2 & 0.67  & 3.4 & 5.2 &  $1.1\times 10^{3}$  \\
\hline
3 & 4 & 0.49 &  1.0 & 1.2 & 58 \\
\hline
4 & 7 & 0.54 & 2.0  & 0.79 & 26 \\
\hline
5 & 10 & 0.38 & 0.93 & 0.48 & 3.5\\
\hline
6 & 14 & 0.25  & 0.66 & 0.25 & 0.64\\
\hline 
7 & 18 & 0.12 & 0.32 & 0.12 & 0.12 \\
\hline
8 & 23 & 0.064 & 0.19 & 0.060  & 0.036 \\
\hline
9 & 28 & 0.032  & 0.093  & 0.032  & 0.012 \\
\hline
10 & 34 & 0.014 & 0.036 & 0.016 & $3.9 \times 10^{-3}$ \\
\hline
\hline
15 & 70 & $1.9\times 10^{-3}$ & $5.1\times 10^{-3}$ & $2.1 \times 10^{-3}$ & $6.8 \times 10^{-5}$ \\
\hline
\end{tabular}

\label{T8}
\end{table}

\subsection{Recovering $a$ and the potential function}

In both the Koiso--Cao and the Wang--Zhu cases, the cohomology parameter $a$ and the coefficient determining the potential dunctions $\phi_{KC}$ and $\phi_{WZ}$ can be considered as variables.  The $T_{1}$ and $T_{2}$ minimisation algorithms were run with $a$ and the coefficient as input variables (the $T_{2}$ algorithm being run over the spaces $\mathcal{S}_{d}^{U(2)}$ in the Koiso--Cao case). In all cases the variable $a$ and the coefficient of the potential function converge to the correct values.  This is a promising sign that the search for quasi-Einstein metrics on $\mathbb{CP}^{2}\sharp 2\overline{\mathbb{CP}}^{2}$ could recover the correct cohomology variable $a$ as well as the parameters determining the conformal factor and the potential function.  
  
\section{Results: Quasi-Einstein metrics} \label{resultsQEM}
\subsection{The L\"u--Page--Pope metric}

The $T_{3}$ minimisation algorithm was performed on the spaces $\mathcal{S}^{U(2)}_{d}$ and $\mathcal{S}^{\mathbb{Z}_{2}}_{d}$.  Here the parameter $m=2$ was used as the resulting metrics yield approximations to a warped product Einstein metric on ${\mathbb{CP}^{2}\sharp \overline{\mathbb{CP}}^{2}\times \mathbb{CP}^{1}}$. As mentioned previously, in this case it is known that
$$a=1, \ b\approx 0.076527, \ c\approx 1.002924, \textrm{ and } d\approx 0.588325.$$ 
As in the case of the Koiso--Cao soliton, the algorithm seems to stall when performed on the spaces $\mathcal{S}^{\mathbb{Z}_{2}}_{d}$.  The results of the algorithm on $\mathcal{S}^{U(2)}_{d}$ are presented in Table \ref{T9}.  Here a similar convergence profile to that of the Koiso--Cao soliton is seen.  Hence there is good evidence that the approximations are really converging to the L\"u--Page--Pope metric with $m=2$.\\
\\
The algorithm was also performed on the spaces $\mathcal{S}^{U(2)}_{d}$ with the parameters $a,b,c$ and $d$ taken to be variables.  In this case the procedure returns approximations where the parameters converge to the values given previously.

\begin{table}[!ht]
\centering
\caption{$T_{3}$ minimisation on the spaces $\mathcal{S}^{U(2)}_{d}$ for the L\"u--Page--Pope metric with $m=2$}
\begin{tabular}{|c|c||c|c||c|c|}
\hline 
Degree $d$ & $N_{d}$ & $\mathcal{E}(T_{4})$  & $Max(T_{4})$ & $\mathcal{E}(T_{3})$ & $Max(T_{3})$  \\
\hline
\hline
2 & 1 & 0.28 & 1.7 &  4.9  & 190  \\
\hline
3 & 2 & 0.22 & 0.40  & 0.53 & 2.9 \\
\hline
4 & 3 & 0.074 & 0.15  & 0.20 & 0.71 \\
\hline
5 & 4 & 0.035  &  0.063 & 0.067 & 0.081 \\
\hline
6 & 5 & 0.012  & 0.024  & 0.023 & 0.011\\
\hline 
7 & 6 & $4.1\times 10^{-3}$ & $8.1\times 10^{-3}$ & $7.5 \times 10^{-3}$ & $1.2\times 10 ^{-3}$ \\
\hline
8 & 7 & $1.3\times 10^{-3}$ & $2.6\times 10^{-3}$ &  $2.4\times 10^{-3}$ &$1.2 \times 10^{-4}$  \\
\hline
9 & 8 & $4.1\times 10^{-4}$ & $7.8\times 10^{-4}$  & $7.7\times 10^{-4}$  & $1.4\times 10^{-5}$ \\
\hline
10 & 9 & $2.4 \times 10^{-4}$ & $5.0\times 10^{-4}$ & $4.4 \times 10^{-4}$ & $4.2\times 10^{-6}$ \\
\hline
\end{tabular}

\label{T9}

\end{table}

\subsection{QEMs on $\mathbb{CP}^{2}\sharp2\overline{\mathbb{CP}}^{2}$}

The $T_{3}$-minimisation algorithm was performed on the spaces $\mathcal{S}^{\mathbb{Z}_{2}}_{d}$. To begin with the parameter $m=2$ and the coefficients $a,b,c$ and $d$ were fixed to be
\begin{equation}
a=2, \ \ b = -0.0744357, \ \ c=1.00482 \textrm{ and } d  = -0.463585. 
\end{equation}
The results of the minimisation algorithm are presented in Table \ref{T10}.  Here there is a similar convergence of the error over the spaces $\mathcal{S}^{\mathbb{Z}_{2}}_{d}$ where $d\leq 9$ to that which occurs for the Wang--Zhu soliton in Table \ref{T8}. However, the convergence slows down; there is no improvement going from $d=10$ and $d=15$. The $d=15$ coefficients were then taken as a starting point for the $T_{3}$ minimisation algorithm with the parameters $a$, $b$ $c$ and $d$ taken as variables. The approximation found in this way  is then taken as the intial condition for higher values of $m$.  The results of this are presented in Table \ref{T11}. The resulting metrics approximately solve the quasi-Einstein equation to the same extent that the $d=15$ approximations to the Wang--Zhu soliton solve the Ricci soliton equation.  It seems that there is evidence that, if conformally K\"ahler quasi-Einstein metrics do exist in the form being investigated, then the cohomology class is not necessarily the canonical one. The final column gives an idea of how well the approximate quasi-Einstein metrics solve the Ricci soliton equation.  One can see that there is evidence that the metrics are converging to the Wang-Zhu soliton as $m\rightarrow \infty$. 
\begin{table}[!ht]
\centering
\caption{$T_{3}$ minimisation on the spaces $\mathcal{S}^{\mathbb{Z}_{2}}_{d}$ with $m=2$ for the manifold $\mathbb{CP}^{2}\sharp 2\overline{\mathbb{CP}}^{2}$}
\begin{tabular}{|c|c||c|c||c|c|}
\hline 
Degree $d$ & $N_{d}$ & $\mathcal{E}(T_{4})$  & $Max(T_{4})$ & $\mathcal{E}(T_{3})$ & $Max(T_{3})$  \\
\hline
\hline
2 & 2 & 0.67 & 3.4 & 5.2  & $1.1\times 10^{3}$ \\
\hline
3 & 4 & 0.49 & 1.1  & 1.2 & 59  \\
\hline
4 & 7 & 0.54  & 2.0  & 0.80 & 26 \\
\hline
5 & 10 & 0.38  & 0.99  & 0.48 & 3.6 \\
\hline
6 & 14 & 0.23 &  0.65 & 0.24 & 0.67\\
\hline 
7 & 18 & 0.12 & 0.29 & 0.12 & 0.11 \\
\hline
8 & 23 & 0.068 & 0.21 &  0.061 & 0.036 \\
\hline
9 & 28 & 0.043 & 0.12 & 0.035 & 0.013 \\
\hline
10 & 34 & 0.032 & 0.063 & 0.024 & $4.7 \times 10^{-3}$  \\
\hline
\hline
15 & 70 & 0.027  & 0.053 & 0.019 & $1.3 \times 10^{-3}$\\
\hline
\end{tabular}
\label{T10}
\end{table}

\begin{table}[!ht]
\centering
\caption{Errors for $d=15$ approximations with no fixed cohomology class}
\begin{tabular}{|c|c|c|c|c|c|c|}
\hline 
$m$ & $a$ & $b$ & $c$ & $d$ & $\mathcal{E}(T_{4})$ & $\mathcal{E}(T_{2})$ \\
\hline
\hline
2 & 1.99256 & -0.0773635 & 0.999924 & -0.445249 & $2.2\times10^{-3}$ & 0.21\\
\hline
3 & 1.99557 & -0.0581581 & 0.999901 & -0.450592 & $1.9\times 10^{-3}$ & 0.13\\
\hline
5 & 1.99771 & -0.0387679 & 0.999819 & -0.456966 & $2.4 \times 10^{-3}$  & 0.065\\
\hline
10 & 1.99897 & -0.0201648 & 0.999956 & -0.488660 & $2.3\times 10^{-3}$ & 0.024\\
\hline
\end{tabular}
\label{T11}
\end{table}
\newpage
\section{Conclusions and Future Work}

We have demonstrated a very straightforward algorithm for approximating the Wang--Zhu soliton on any Fano toric K\"ahler manifold.  As discussed, there is some numerical evidence that the algorithm based on minimising the Ricci curvature quantity $T_{2}$ can stall.  Hence a more sophisticated method of optimisation could be used.  One method might be to use the polytope $P$ to define a number of charts on the manifold where the Ricci tensor (and other tensor quantities) does not become singular.  The integrals in the algorithm could then be computed in each chart.  Another method would be to use a more sophisticated minimisation algorithm in conjunction with the Levenberg-Marquardt method to prevent the algorithm being trapped in local minima.\\
\\
The evidence for quasi-Einstein metrics on $\mathbb{CP}^{2}\sharp 2\overline{\mathbb{CP}}^{2}$ is intriguing.  It seems that the algorithm does not find a straightforward generalisation of the L\"u--Page--Pope metric where the K\"ahler class is a representative of $c_{1}$. This evidence is strengthened by the fact that much better convergence is found (Table \ref{T11}) if the cohomology class is allowed to vary. However, given the evidence that the Ricci curvature minimisation methods can stall, the failure of the method to converge might be algorithmic. There is also the possibility that there are quasi-Einstein metrics not of the form considered (conformally K\"ahler with $J$-invariant Ricci tensor). However, searches using minor modifications of the algorithms to search for metrics without $J$-invariant Ricci tensor did not yield improvements. The main question arising from this work would be to give a proof of existence of quasi-Einstein metrics on $\mathbb{CP}^{2}\sharp 2\overline{\mathbb{CP}}^{2}$ or to find an obstruction. 

\bibliography{QEMRefs}{}
\bibliographystyle{acm}


\end{document}